\def\versiondate{3 Sep. 2016}
\input math.macros

\let\nobibtex = t
\let\noarrow = t
\input eplain.tex
\beginpackages
\usepackage{url}
\endpackages

\input Ref.macros
\input lrlEPSfig.macros

\checkdefinedreferencetrue
\continuousfigurenumberingtrue
\theoremcountingtrue
\sectionnumberstrue
\forwardreferencetrue
\citationgenerationtrue
\nobracketcittrue
\hyperstrue
\initialeqmacro

\input\jobname.key
\bibsty{myapalike}

\def\verts{{\ss V}}
\def\vertex{{\ss V}}
\def\Td{{\Bbb T}_d}  
\def\Tdo{{\bf o}}  
\def\T{{\Bbb T}}  
\def\Fr{{\Bbb F}}  
\def\F{{\scr F}}  
\def\gp{\Gamma}
\def\gpe{\gamma}
\def\gh{G}

\def\edges{{\ss E}}
\def\bp{o}

\def\ed#1{[#1]}  
\def\uedg(#1, #2){\ed{#1, #2}}
\def\oedg(#1){\Seq{#1}}  
\def\leb{{\cal L}}  
\def\Gm{{\frak m}}  
\def\dbar{\bar d}
\def\dtbar{\bar d_2}
\def\Aut{{\rm Aut}}
\def\pf{\mu^{\rm pm}}  
\def\dcol{\mu^{\rm col}}  
\def\mc{\mu^{\rm mc}}
\def\sphsum{{\bf \Sigma}}
\def\sgn{\mathop{\rm sgn}}
\def\Corr{\mathop{\rm Corr}}
\def\SD{\mathop{\rm SD}}
\def\sz{\mathop{\ss size}}
\def\minbi{\mathop{\ss MinBi}}
\def\maxbi{\mathop{\ss MaxBi}}
\def\minbifac{\mathop{\ss MinBiFac}}
\def\maxbifac{\mathop{\ss MaxBiFac}}
\def\wf{\mathop{\ss WF}}
\def\cut{\Pi}

\def\eqD{\buildrel {\cal D} \over =}  
\def\P{{\Bbb P}}

\def\A{{\frak A}}  
\def\parisi{{\ss P_{\! *}}}

\ifproofmode \relax \else\head{To appear in {\it Combin. Probab.
Comput.}}
{Version of \versiondate}\fi 
\vglue20pt

\title{Factors of IID on Trees}

\author{Russell Lyons}

\abstract{Classical ergodic theory for integer-group actions uses entropy
as a complete invariant for isomorphism of IID (independent, identically
distributed) processes (a.k.a.\ product measures). This theory holds for
amenable groups as well.  Despite recent spectacular progress of Bowen, the
situation for non-amenable groups, including free groups, is still largely
mysterious. We present some illustrative results and open questions on free
groups, which are particularly interesting in combinatorics, statistical
physics, and probability.  Our results include bounds on minimum and maximum
bisection for random cubic graphs that improve on all past bounds.
}

\bottomIII{Primary 
05C30, 
05C70, 
05C80, 
37A35, 
37A50, 
60G10. 
Secondary
60G15. 
}
{Non-amenable groups, regular graphs, bisection, local rules, Ising, tail
field.}
{Research partially supported by Microsoft Research and
NSF grant DMS-1007244.}

\bsection{Introduction}{s.intro}

Let $\gp$ be a group and $X$ and $Y$ be two sets on which $\gp$ acts. A map
$\phi \colon X \to Y$ is called \dfn{$\gp$-equivariant} if $\phi$
intertwines the actions of $\gp$: 
$$
\phi(\gpe x) = \gpe \big(\phi(x)\big) \qquad (\gpe \in \gp,\, x \in X)
\,.
$$
If $X$ and $Y$ are both measurable spaces, then a $\gp$-equivariant
measurable $\phi$ is called a \dfn{$\gp$-factor}.
Let $\mu$ be a measure on $X$. If $\phi$ is a $\gp$-factor, then the
push-forward measure $\phi_* \mu$ is called a \dfn{$\gp$-factor of $\mu$}.
The measure $\mu$ is \dfn{$\gp$-invariant} if 
$$
\mu(\gpe B) = \mu(B) \qquad (\gpe \in \gp,\, B \subseteq X \hbox{
measurable}).
$$
If $\nu$ is a measure on $Y$, then a $\nu$-a.e.-invertible $\gp$-factor
$\phi$ such that $\nu = \phi_* \mu$ is called an \dfn{isomorphism} from
$(X, \mu, \gp)$ to $(Y, \nu, \gp)$.
Classical ergodic theory is concerned with the case that $\gp = \Z$.
In probability theory, we often have that $X$ and $Y$ are product spaces
of the form $\A^\gp$ or, more generally, $\A^W$, where $\A$ is a measurable
space, called the \dfn{base}, and $W$ is a countable set on which $\gp$ acts.
Note that in this case, $\gp$ acts on $\A^W$ by 
$$
\big(\gpe \omega\big)(x) 
:=
\omega( \gpe^{-1} x)  \qquad (\omega \in \A^W,\, x \in  W,\, \gpe
\in \gp)
\,.
$$

When $\A$ is finite, $\lambda$ is the uniform measure on $\A$, and $\gp$ is
denumerable, then $(\A^\gp,
\lambda^\gp, \gp)$ is called the \dfn{$|\A|$-shift over $\gp$}.
When $\mu$ and $\nu$ are more general
product measures with finite base spaces and
$\gp = \Z$, the problem of whether $(\A^\gp, \mu, \gp)$ and $(\A^\gp, \nu,
\gp)$ are isomorphic is very
old. It was solved through the introduction of entropy by
\refbmulti{Kolmogorov:new,Kolmogorov:entropy} and \ref b.Sinai:ent/, and
the work of \ref b.Orn:same/.
In particular, the entropy of a $k$-shift is $\log k$. 
Factors play a key role in this and other aspects of ergodic theory.

The theory of entropy and its applications
was extended to amenable groups by \ref b.OrnW:amen/.
One important feature is that the entropy of a factor of an invariant
probability measure $\mu$ is at most the entropy of $\mu$.
\ref b.OrnW:amen/ noted that no reasonable definition of entropy on free
groups of rank at least 2 has this property, since, as they showed, the
4-shift is a factor of the 2-shift over such groups.
Clearly the 2-shift is also a factor of the 4-shift.
The problem whether the 2-shift and the 4-shift, for example, are
isomorphic was finally solved by \ref b.Bowen:fg/, who introduced a notion
of entropy for free-group actions that is invariant under isomorphism.
His notion of entropy again assigns the value
$\log k$ to a $k$-shift over any free group.
Bowen's work is the analogue of that of Kolmogorov and Sina\u\i.
However, most of Ornstein theory remains {\it terra incognita}.

One reason, therefore, to study factors over free groups, and especially
factors of product measures, is to understand how to extend Ornstein
theory. 
Other reasons arise from questions in probability theory, combinatorics,
and computer science, as well as the ergodic theory of equivalence
relations.
We shall discuss questions from most of these areas here, providing some
results and highlighting some particularly interesting open questions.

Recent papers concerning factors (generally of IID processes or their
continuous analogue, Poisson point processes)
that involve a mix of probability and
combinatorics include
\ref b.HP:match/,
\ref b.Timar:pp/, 
\ref b.Ball:thin/,
\ref b.HPS:match/,
\ref b.CPPR:grav/,
\ref b.Soo:match/, \ref b.Hol:geommatch/,
\ref b.HLS:split/, 
\ref b.Timar:vor/, 
\ref b.Mester:uniform/,
\ref b.LN:factor/,
\ref b.ABGMP:map/,
\ref b.CL:match/,
\ref b.GamSudan/,
\ref b.CGHV/,
\ref b.HarVir:transtv/,
\ref b.BSV:decay/,
\ref b.Conley:Brooks/,
\ref b.Kun:Lip/,
\ref b.QS:monotone/,
and
\ref b.GGP:thick/.

The utility of factors of IID processes on non-amenable groups 
has
been shown in various ways. For example, see 
\ref b.Popa/,
\ref b.ChifanIoana/,
\ref b.Houdayer/,
\ref b.Lyons:fixed/,
\ref b.AW:Bernoulli/,
and
\ref b.Kun:Lip/.

Because in many situations one has the natural Cayley graph of a free
group, that is, a regular tree, one also is often interested in processes
that are invariant under the full automorphism group of the tree and,
similarly, in factors that are equivariant with respect to the full
automorphism group.
Thus,
let $\Td$ be a $d$-regular tree with $d \ge 3$.
When $d$ is even, this is a Cayley graph of the free group $\Fr_{d/2}$ on
$d/2$ generators. In all cases, it is a Cayley graph of the free product of
$d$ copies of $\Z_2$.

Our greatest interest in this paper is $\Aut(\Td)$-factors 
$$
\phi \colon \big([0, 1]^{\verts(\Td)}, \leb^{\verts(\Td)}\big) \to
\{0, 1\}^{\verts(\Td)}
\,,
$$
where $\leb$ is Lebesgue measure on $[0, 1]$.
We shall often leave off the prefix $\Aut(\Td)$ from the word ``factor".
Since the domain space is product measure, or IID,
such a $\phi$ is called a \dfn{factor of IID}, or \dfn{FIID} for short.
The elements of the domain space are sometimes called \dfn{labels}.
The push-forward measure $\phi_* \leb^{\verts(\Td)}$ is also called an
\dfn{FIID}.
Under the same rubric we shall consider other product measures over either
$\verts(\Td)$ or $\edges(\Td)$, with the codomain also being other product
measurable spaces over either $\verts(\Td)$ or $\edges(\Td)$.

Extending the fundamental example of \ref b.OrnW:amen/, \ref b.Ball:factor/
showed that the 4-shift is an $\Aut(\Td)$-factor of the 2-shift, as is
$\leb^{\verts(\Td)}$.
For this reason, it matters little which product measure is used as the
domain of a factor.

Our contributions in this area are to exhibit weak* limits of FIID
processes that are not themselves FIID; and to use FIID processes in order
to improve on existing bounds for minimum and maximum bisection of random
regular graphs.
For example, \ref b.MonPre/ showed that random 3-regular graphs
asymptotically have bisection width at most 1/6, which we improve to
0.1623.

\bsection{Factors on Trees}{s.trees}

Let $\Tdo$ denote a fixed vertex, the \dfn{root}, of $\Td$.
There is a correspondence between $\Aut(\Td)$-factors of
IID, $\phi$, and spherically symmetric measurable functions
$F \colon [0, 1]^{\verts(\Td)} \to \{0, 1\}$, namely, 
$$
F(\omega) = \big(\phi(\omega)\big)(\Tdo) \qquad (\omega \in [0, 1]^\verts)
\label e.defF
$$
in one direction and 
$$
\big(\phi(\omega)\big)(\gpe^{-1} \Tdo) 
=
F(\gpe \omega) \qquad (\omega \in [0, 1]^\verts,\, \gpe \in \Aut(\Td))
\label e.defphi
$$
in the other.

For a measurable space $\A$, write $\pi_x \colon \A^\verts \to \A$ for the
natural coordinate projections ($x \in \verts$).
For $K \subseteq \verts$, write $\F(K)$ for the $\sigma$-field on
$\A^\verts$ generated the maps $\pi_x$ for $x \in K$. 
Let $B_r := B_r(\Tdo)$ be the graph induced on the set of vertices within
graph distance $r$ of $\Tdo$.
We may approximate $F$ as in \ref e.defF/ 
by spherically symmetric measurable maps
$F_r \colon [0, 1]^{\verts(B_r(\Tdo))} \to \{0, 1\}$ that converge to $F$
a.s. For example, the conditional expectations $F_r := \Ebig{F \bigm| 
\F\big(B_r(\Tdo)\big)}$ converge to $F$ by L\'evy's 0-1 Law.
These maps $F_r$ determine FIIDs
$\phi_r$ via \ref e.defphi/,
called \dfn{block factors of IID} or \dfn{local rules}. We have
that $\phi_r$ converges to $\phi$ a.s.\ (in the product topology)
and therefore $(\phi_r)_*
\leb^\verts$ converges to $\phi_* \leb^\verts$
in Ornstein's $\dbar$-metric.
This metric is defined as follows.
Let $\mu_1$ and $\mu_2$ be two $\gp$-invariant probability measures on $\A^W$,
where $\gp$ acts quasi-transitively on $W$ and $\A$ is a finite set. Let
$W'$ be a section of $\gp\backslash W$. Then 
$$
\dbar(\mu_1, \mu_2)
:=
\min \Big\{\sum_{w \in W'} \Pbig{X_1(w) \ne X_2(w)} \st X_1 \sim \mu_1,\, X_2
\sim \mu_2,\, (X_1, X_2) \hbox{ is $\gp$-invariant} \Big\}
\,.
$$

The \dfn{tail}
$\sigma$-field is defined to be $\bigcap_r \F(\verts \setminus B_r)$.
For $x \in \verts$, let $D_x$ denote the set of vertices separated from
$\Tdo$ by $x$.
If $(x_1, x_2, \ldots)$ is a simple path of vertices in $\Td$, the
corresponding \dfn{1-ended tail} $\sigma$-field is $\bigcap_n
\F\big(D(x_n)\big)$.
Let $\Aut_+(\Td)$ denote the parity-preserving subgroup of $\Aut(\Td)$, i.e.,
$\Aut_+(\Td) := \big\{\gpe \in \Aut(\Td) \st \all {x \in \vertex(\Td)} d(x,
\gpe x) \in 2\N\big\}$, where $d(x, y)$ is the graph distance between $x$
and $y$.
Every $\Aut(\Td)$-invariant $\Aut(\Td)$-ergodic probability measure on
$\A^\verts$ is an equal mixture of two $\Aut_+(\Td)$-invariant
$\Aut_+(\Td)$-ergodic probability measures, and the latter have
trivial 1-ended tail $\sigma$-fields, as shown by \ref b.Pemantle/.
By virtue of being $\Aut_+(\Td)$-ergodic, every FIID has trivial 1-ended tails.

A probability measure $\mu$ on $\A^\verts$ is called
\dfn{$m$-dependent} if $\F(K_1), \dots, \F(K_p)$ are independent whenever
the sets $K_i$ are pairwise separated by graph distance $> m$. We say that
$\mu$ is \dfn{finitely dependent} if it is $m$-dependent for some $m <
\infty$.
For example, a block FIID that depends on the ball of radius $r$ is
$2r$-dependent.

According to the Kolmogorov 0-1 Law, the tail $\sigma$-field is trivial for
every IID probability measure. Reasoning similar to its proof shows
the second of the following implications for $\Aut(\Td)$-invariant
processes:
$$
\hbox{block FIID} \implies \hbox{finitely dependent} \implies \hbox{trivial
tail}.
$$
It is open whether finitely dependent implies FIID and whether trivial tail
implies FIID.
These questions are resolved on $\Z$:
finitely dependent implies FIID by using the VWB condition of
\ref b.Orn:book/ and trivial tail does not imply FIID (even for finite $\A$)
by \ref b.Orn:KnotB/ and \ref b.Kalikow:TTinv/. It is
also known that for finite $\A$,
FIID implies trivial 1-ended tail $\sigma$-fields, as proved by \ref b.RohSin/.
This latter implication is false for $\A = [0, 1]$.
We note, however, that \ref b.Smoro/ proved that
Gaussian processes on $\Z$ with trivial 1-ended tail {\it are} FIID.

The following question is due to \ref b.Bowen:FIIDiso/:

\procl q.FIIDiso
Is every FIID process isomorphic to an IID process?
\endprocl

\ref b.Orn:factor/ proved this holds on $\Z$.
It does not suffice on $\Td$ to have factor
maps each way, since this holds for the 2-shift and 4-shift, but these are
not, by \ref b.Bowen:fg/, isomorphic. 
Note that \ref b.Popa/ proved that there exist non-amenable groups where
FIIDs are not necessarily isomorphic to IIDs.

Many of the above questions can be asked about invariant processes on
non-amenable groups more generally, not just free groups.

\procl q.findep
Is every finitely dependent process an FIID?
\endprocl

This holds in the amenable case again by using the VWB condition, here
defined by \ref b.Adams:VWB/.

We now present an example of an FIID on $\Td$ with finite $\A$
whose tail $\sigma$-field is full (everything).
Such examples on $\Z$ were given by \ref b.OrnW:bilat/ (who proved that
every process is isomorphic to one whose tail $\sigma$-field is full), \ref
b.BDS:bilat/, and \ref b.BurSt:surf/.

\procl p.pm 
There exists a unique $\Aut(\Td)$-invariant probability measure, $\pf$, on
the set of perfect matchings of $\Td$; it is an FIID whose tail
$\sigma$-field is full.
\endprocl

\proof
Since the stabilizer $\gp$ of $\Tdo$ in $\Aut(\Td)$ acts transitively 
on the set of perfect matchings of $\Td$,
there is a unique $\gp$-invariant probability measure, $\pf$,
on the set of perfect matchings.
This measure is easy to construct by
starting at $\Tdo$, choosing uniformly 
at random one of its $d$ incident edges to be in the
matching, and then working outwards independently, where every time there
is a choice between $d-1$ edges, they are equally likely to be in the
matching.
Using the independence, it is not hard to see that $\pf$ is actually
$\Aut(\Td)$-invariant.
Although it is far from obvious, $\pf$ is an FIID, as shown by \ref
b.LN:factor/.

To see that the tail is full, consider any event $A$ of perfect matchings
and any radius $r \ge
0$. Let $A_r$ be the event consisting of all perfect matchings that agree
with some element of $A$ when restricted to the complement of
$\edges(B_r(\Tdo))$. We claim that $A
= A_r$ for all $r$, which will imply that the tail of $\pf$ is full.
We prove this by induction on $r$. It is clear that $A_0 = A$.
Now let $\omega \in A$ and $r \ge 0$. By definition, there exists some
$\omega' \in A_{r+1}$ that agrees with $\omega$ outside $B_{r+1}(\Tdo)$.
Consider an edge $e \in B_{r+1}(\Tdo) \setminus B_{r}(\Tdo)$.
Let $F$ be the set of $d-1$ edges incident to $e$ that do not lie in
$B_{r+1}(\Tdo)$.
Since $\omega(e) = 1$ iff $\omega(f) = 0$ for all $f \in F$, and likewise
for $\omega'$, it follows that $\omega(e) = \omega'(e)$, whence that
$\omega$ agrees with $\omega'$ outside
$B_r(\Tdo)$, i.e., that $\omega' \in A_r$. Therefore, $A = A_r$ implies that $A =
A_{r+1}$, which completes the induction.
\Qed

For similar reasons, there is a unique $\Aut(\Td)$-invariant probability
measure, $\dcol$, on the set of proper $d$-colorings
of $\edges(\Td)$. This measure is again easy
to construct by
working outwards from $\Tdo$.
Proper $d$-colorings can also be regarded as Cayley diagrams of the free
product, $\Z_2^{*d}$, of $d$ copies of $\Z_2$.

\procl q.Cayleycolor
Is $\dcol$ an FIID? 
\endprocl

This is open. A positive answer would imply that the set of
$\Aut(\Td)$-factors is equal to the set of $\Aut(\Td)$-invariant
$\Z_2^{*d}$-factors. Note that every $\Z_2^{*d}$-invariant probability measure
induces an $\Aut(\Td)$-invariant probability measure by averaging with
respect to the stabilizer of $\Tdo$ in $\Aut(\Td)$.

Let $1, \ldots, d$ be the $d$ colors we use.
It is also open whether $\dcol = \phi_* \pf$ for some $\Aut(\Td)$-factor
$\phi$ that colors every edge in the perfect matching by color $1$, i.e.,
$\big(\phi(\omega)\big)(e) = 1$ for all $e$ with $\omega(e) = 1$.

We mention a partial result towards answering \ref q.Cayleycolor/.
A proper $d$-coloring is the same as a list $(P_1, \ldots, P_d)$ of $d$
disjoint perfect matchings.
It is possible to obtain as an FIID a probability measure on lists $(P_1,
\ldots, P_{d-2}, \{Q_1, Q_2\})$, where $P_i$ and $Q_j$ are disjoint perfect
matchings, but $\{Q_1, Q_2\}$ is unordered.
Indeed, choose $P_1$ via the FIID $\pf$.
Note that $\big([0, 1], \leb\big)$ is isomorphic to $\big([0, 1]^{\N},
\leb^{\N}\big)$, so that when we create $P_1$ as an FIID, we may use only
the first coordinates of the labels, reserving the later coordinates for
further use.
Deleting the edges of $P_1$ decomposes $\Td$ into a forest of copies of
$\T_{d-1}$. Provided $d-1 \ge 3$, we may choose perfect matchings in each
copy by using the second coordinates of the labels and let $P_2$ be their
union.
This procedure may be continued until we are left with trees of degree 2.
Each such tree is decomposed uniquely as a set of 2 perfect matchings.
We must decide, given a perfect matching $P$ of a tree $T$ of degree 2 and
a perfect matching $P'$ of another tree $T'$ of degree 2, whether $P$ and
$P'$ belong to the same $Q_i$ or not. In order to make this decision for
all perfect matchings and all trees, it suffices to make the decision for
pairs of trees that are at distance 1 from each other in $\Td$. In such a
case, there is a unique edge $e$ that is incident to both trees. Let $e_1$
and $e_2$ be the two edges in $T$ that are adjacent to $e$ and $e'_1$ and
$e'_2$ be the two edges in $T'$ that are adjacent to $e$. Let $U_i$ and
$U'_i$ be the corresponding labels of these edges ($i = 1, 2$).
Then let the perfect
matching containing $e_1$ belong to the same $Q_i$ as the perfect matching
containing $e'_1$ iff $(U_1-U_2)(U'_1-U'_2) > 0$.

\bsection{Tree-Indexed Markov Chains and Ising Measures}{s.ising}

Next we consider the simplest types of invariant processes after IID,
namely, 2-state symmetric $\Td$-indexed Markov chains.
Let $|\theta| \le 1$ and consider the transition matrix 
$$
\left(\matrix{{1+\theta\over2}&{1-\theta\over2}\cr
               {1-\theta\over2}&{1+\theta\over2}\cr}\right).
$$
For $\theta \ge 0$, another way to think of this transition matrix, which
explains this parametrization, is to keep the same state with probability
$\theta$ and to choose uniformly among the two states independently of the
current state with probability $1 - \theta$.
For $\theta \le 0$, the interpretation is slightly different:
change to the opposite state with probability
$|\theta|$ and to choose uniformly among the two states independently of the
current state with probability $1 - |\theta|$.
The tree-indexed Markov chain $\mc_\theta$ is obtained by assigning to the
root one of the 2 states with equal probability, then proceeding to the
neighbors of the root by using an independent transition from the above
matrix, etc.
When the two states are $\pm 1$, this is known as the free Ising measure on
$\Td$, ferromagnetic when $\theta \ge 0$. In this case, the states are
known as spins. We shall use this terminology for convenience.

The description of $\mc_\theta$ does not make it apparent that $\mc_\theta$
is an invariant measure, but it is not hard to check that it is indeed
invariant.  However, an important alternative description makes this
invariance obvious. Namely, consider the clusters of Bernoulli($|\theta|$) bond
percolation on $\Td$.  If $\theta \ge 0$, then for each cluster, assign all
vertices the same spin, with probability 1/2 for each spin,
independently for different clusters. If
$\theta \le 0$, then assign each cluster one of its two proper $\pm
1$-colorings, with probability 1/2 each, independently for different
clusters. It is easy to see that this gives $\mc_\theta$.

It is known that $\mc_\theta$ has a trivial tail iff $|\theta| \le 1/\sqrt
{d-1}$.
It is also known that $\mc_\theta$ is an FIID if $|\theta| \le 1/(d-1)$,
but is not an FIID if $|\theta| > 1/\sqrt{d-1}$.
It is open whether $\mc_\theta$ is an FIID for $1/(d-1) < |\theta| \le
1/\sqrt{d-1}$.
It is also open whether there is a critical $\theta_0$ such that
$\mc_\theta$ is an FIID for $0 \le \theta < \theta_0$ and not an FIID for
$\theta_0 < \theta \le 1$; the analogous question is also open for $\theta
< 0$.
The history of the result for tail triviality is reviewed in Sec.~2.2 of
\ref b.EKPS/.

The fact that $\mc_\theta$ is an FIID for $|\theta| \le 1/(d-1)$ is 
easy to see: In this regime, all $|\theta|$-clusters are finite
a.s. Let $U(e)$ and $U_i(x)$
be IID uniform $[0, 1]$ random variables for $e \in \edges(\Td)$, $i \in
\{1, 2\}$, and $x \in \verts(\Td)$.
Choose the $|\theta|$-clusters by using the edges with $U(e) \le |\theta|$.
Given a cluster $C$, let its vertex with the minimum $U_1(x)$ be $x_C$ and let
the spins in $C$ equal $\sgn \big(U_2(x_C) - 1/2\big)$ if $\theta \ge 0$,
while if $\theta < 0$, let the spins in $C$ equal the proper $\pm
1$-coloring whose spin at $x_C$ equals $\sgn \big(U_2(x_C) - 1/2\big)$.

In unpublished work, this author and later Lewis Bowen gave values
$\theta_d$ such that for $|\theta| > \theta_d$, the measure $\mc_d$ is not
an FIID. This was improved to $\theta_d = 1/\sqrt{d-1}$ by \ref
b.Sly:FIID/, but his proof was not published. We give that proof here
because we shall adapt it to prove other results as well. This value of
$\theta_d$ can also be established by using a result of \ref b.BSV:decay/,
which characterizes the rate of decay of the correlation of $\sigma(\Tdo)$
and $\sigma(x)$ as the distance between $\Tdo$ and $x$ tends to infinity,
where $\sigma$ is any FIID whose values at the vertices are real valued and
square integrable.
In particular, the correlation is at most $\big(n(d-2)/d +
1\big)/(d-1)^{n/2}$ in absolute value when the distance is $n$.
Of course, this holds as well for weak* limits of FIID processes.
In particular, weak* limits of FIID are strongly mixing, while on $\Z$,
they need not even be ergodic.
Note that $\mc_\theta$ has trivial 1-ended tails for all $|\theta| < 1$.

The following is at the heart of Sly's proof, with the last observation
about $\dtbar$-closure due to this author and Peres in 2013.
Here, given two invariant probability measures $\mu_1$ and $\mu_2$ on
$\R^{\verts}$, we define 
$$
\dtbar(\mu_1, \mu_2)
:=
\min \Big\{ \Ebig{|X_1(\bp)-X_2(\bp)|^2}^{1/2} \st X_1 \sim \mu_1,\, X_2 \sim
\mu_2,\ (X_1, X_2) \hbox{ is $\gp$-invariant} \Big\}
\,.
$$
Note that FIID processes whose 1-dimensional marginals have
finite second moments
are $\dtbar$-limits of block factors.

\procl t.sly Let $\gh$ be a graph for which there is some
unimodular group $\gp$ of
automorphisms that acts transitively on $\verts(\gh)$.
Let $\bp \in \verts(\gh)$. Write $S_n$ for the set of vertices at distance
$n$ from $\bp$.
Suppose that $x \mapsto \sigma(x)$ ($x \in \verts(\gh)$)
is a $\gp$-invariant process with law $\mu$
on $\R^{\verts(\gh)}$. Assume that 
$0 < \Var\big(\sigma(\bp)\big) < \infty$.
Define $\sphsum_n := \sum_{x \in S_n} \sigma(x)$.
If $\lim_{n \to\infty} \Var(\sphsum_n)/|S_n| = \infty$ and
$\limsup_{n \to\infty} \big|\Corr\big(\sigma(\bp), \sphsum_n\big)\big| >
0$, then $\mu$ is not a $\gp$-equivariant FIID, nor is $\mu$ in the
$\dtbar$-closure of the finitely dependent processes.
\endprocl

Note that the condition $\limsup_{n \to\infty} \big|\Corr\big(\sigma(\bp),
\sphsum_n\big)\big| > 0$ alone implies that $\mu$ has a non-trivial tail.

\proof
Without loss of generality, we may assume that $\Ebig{\sigma(\bp)} = 0$ and
that $\SD\big(\sigma(\bp)\big) = 1$.
We shall show that if
$\lim_{n \to\infty} \Var(\sphsum_n)/|S_n| = \infty$ and
$\mu$ lies in the
$\dtbar$-closure of the finitely dependent processes, then
$\lim_{n \to\infty} \Corr\big(\sigma(\bp), \sphsum_n\big) = 0$.

Let $\epsilon > 0$.
Choose an invariant process $(X, Y)$ on $\R^\verts \times \R^\verts$
such that $X \sim \mu$, $Y$ is finitely dependent,
and $\SD\big(X(\bp)-Y(\bp)\big) < \epsilon$.
For simplicity of notation, we take $X = \sigma$.
Write $\sphsum_n^Y := \sum_{x \in S_n} Y(x)$.
By the Mass-Transport Principle, we have that $\Ebig{Y(\bp) \sphsum_n} =
\Ebig{\sigma(\bp) \sphsum_n^Y}$.
By finite dependence, we have $\Var\big(\sphsum_n^Y\big) = O\big(|S_n|\big) =
o\big(\Var(\sphsum_n)\big)$, whence 
$$
\EBig{Y(\bp) \sphsum_n/\SD\big(\sphsum_n\big)} =
\EBig{\sigma(\bp) \sphsum_n^Y/\SD\big(\sphsum_n\big)} \to 0
$$
as $n \to\infty$.
By the Cauchy--Schwarz inequality, we have
$$
\Big|\EBig{Y(\bp) \sphsum_n/\SD\big(\sphsum_n\big)} - 
\EBig{\sigma(\bp) \sphsum_n/\SD\big(\sphsum_n\big)}\Big|
\le
\SD\big(Y(\bp)-\sigma(\bp)\big) < \epsilon
\,.
$$
Taking $n \to\infty$, it follows that
$$
\limsup_{n \to\infty} \big|\Corr\big(\sigma(\bp), \sphsum_n\big)\big|
=
\limsup_{n \to\infty}
\Big|\EBig{\sigma(\bp) \sphsum_n/\SD\big(\sphsum_n\big)}\Big|
<
\epsilon
\,,
$$
as desired.
\Qed

The following is Sly's result.

\procl c.MCFIID
For $|\theta| > 1/\sqrt{d-1}$, the $\Td$-indexed Markov chain $\mc_\theta$
is not an FIID.
\endprocl

\proof
Let $|\theta| > 1/\sqrt{d-1}$
and $\sigma \sim \mc_\theta$.
We verify the conditions of \ref t.sly/.
Note that 
$$
\Corr\big(\sigma(x), \sigma(y)\big) =
\Ebig{\sigma(x) \sigma(y)} = \theta^n
$$
when $x$ and $y$ are at distance $n$ from each other.
Also, $|S_n| = d (d-1)^{n-1}$.
Therefore, $\Ebig{\sigma(\Tdo) \sphsum_n} = |S_n| \theta^n$ and 
$$
\Var(\sphsum_n)/|S_n| 
=
1 + \sum_{k=1}^{n-1} (d-2)(d-1)^{k-1} \theta^{2k} + (d-1)^n \theta^{2n}
\sim c |S_n| \theta^{2n}
$$
for some constant $c$ as $n \to\infty$.
Hence the conditions of \ref t.sly/ follow.
\Qed

It is open whether discrete FIID processes are closed in the
$\dbar$-topology, as they are on $\Z$ (see \ref b.Orn:book/).
We shall use \ref t.sly/ to show that the class of FIID processes is not
closed in the weak* topology, as is easy to show on $\Z$. 
This was also shown independently by \ref b.HarVir:transtv/ on all infinite
finitely generated groups, but their proof does not show the same for
discrete processes. %

Define $\rho_d := 2\sqrt{d-1}/d$. The spectrum of the transition operator
for simple random walk on $\Td$ is the interval $[-\rho_d, \rho_d]$, as
shown by \ref b.Kesten:amenB/. 
Let $\sigma$ be the Gaussian wave function of \ref b.CGHV/ 
with eigenvalue $\rho$ for the transition operator; this is a centered
Gaussian field on $\Td$ whose covariances satisfy the recurrence 
$$
c_0 = 1;\
c_1 = \rho;\
(d-1) c_{k+1} - d \rho c_k + c_{k-1} = 0 \ (k \ge 1)
\,,
$$
where $c_k$ is the covariance between each pair of vertices at distance $k$.

\procl c.w*FIID
There is a Gaussian process on $\Td$ that is not an FIID but is a weak*
limit of FIID processes. There is a $\{0, 1\}$-valued process that is not
an FIID but is a weak* limit of FIID processes. 
\endprocl

\proof
Let $\sigma$ be the Gaussian wave function 
with eigenvalue $\rho_d$ for the transition operator; \ref b.CGHV/ show
that $\sigma$ is a weak* limit of FIIDs.
Induction shows that
$$
\Corr\big(\sigma(x), \sigma(y)\big) = \big(n(d-2)/d + 1\big) (d-1)^{-n/2}
$$
for $x$ and $y$ at distance $n$.
Therefore, $\Ebig{\sigma(\Tdo) \sphsum_n} \sim c n |S_n|^{1/2}$ and 
$
\Var(\sphsum_n)/|S_n| 
\sim
c' n^2
$
for some positive constants $c$ and $c'$ as $n \to\infty$.
Hence the conditions of \ref t.sly/ follow.

Now let $\tau := \sgn \sigma$.
Since $\sigma$ is a weak* limit of FIIDs, so is $\tau$.
Note that there exist positive constants $c_1$ and $c_2$ such that
if $Z_1$ and $Z_2$ are jointly normal random variables, then
$$
c_1 |\Corr(\sgn Z_1, \sgn Z_2)| \le
|\Corr(Z_1, Z_2)| \le c_2 |\Corr(\sgn Z_1, \sgn Z_2)|
\,.
$$
Hence the above calculations for $\sigma$ hold (up to bounded factors) for
$\tau$ as well.
\Qed

\bsection{Edge Cuts in Finite Graphs}{s.cuts}

Weak* limits of FIID processes on $\Td$ can be used to bound combinatorial
quantities on random $d$-regular graphs or on $d$-regular graphs whose girth
tends to infinity. More generally, they can be used on finite graphs whose
random weak limit is $\Td$. To explain this widely known idea, we first
define ``random weak limit" (for this restricted case).

For a vertex $x$ in a graph $G$, let $B_r(x; G)$ denote the subgraph
induced by the vertices in $G$ whose distance from $x$ is at most $r$.
We consider this subgraph as rooted at $x$.

Let $\Seq{G_n}$ be a sequence of finite graphs. For each $r \ge 1$, let
$p_{n, r}$ denote the probability that a uniformly random vertex $x$ in $G_n$ 
satisfies the property that $B_r(x; G_n)$ is rooted isomorphic to
$B_r(\Tdo; \Td)$, i.e., there is a graph isomorphism from $B_r(x; G_n)$ to
$B_r(\Tdo; \Td)$ that sends $x$ to $\Tdo$.
We say that the \dfn{random weak limit} of $\Seq{G_n}$ is $\Td$ if $\lim_{n
\to\infty} p_{n, r} = 1$ for every $r\ge 1$.
It is evident that every sequence of $d$-regular graphs whose girth tends
to infinity has this property. It is well known that if $G_n$ is a
uniformly random $d$-regular graph on $n$ vertices (or, if $d$ is odd, on
$2n$ vertices), then also $\Seq{G_n}$ has this property with probability 1.
Other terms for this same concept are ``Benjamini--Schramm convergence" and
``local weak convergence".

Now, for the sake of concreteness,
suppose that $\phi$ is a block FIID on $\Td$ associated to the
spherically symmetric measurable map 
$F \colon [0, 1]^{\verts(B_r(\Tdo; \Td))} \to \{0, 1\}$.
Given a graph $G$, one may assign independent uniform $[0, 1]$ random
variables to its vertices and then apply $F$ at every vertex $x$ for
which $B_r(x; G)$ is rooted isomorphic to $B_r(\Tdo; \Td)$.
At other vertices, assign the value 0.
In this way, we obtain a probability measure on $\{0, 1\}^{\verts(G)}$ that
is ``close" to $\phi_* \leb^{\verts(\Td)}$ when $G$ is ``close" to $\Td$.
In particular, the expected number of vertices assigned the value 1 will be
close to $\Pbig{\phi(\cbuldot)(\Tdo) = 1}$.
Informally, we say that $\phi$ is \dfn{emulated} on $G$.

If we want to bound the number of vertices assigned 1 under some constraint
on the set assigned 1, then exhibiting a random set obtained by emulating a
block factor will help. Moreover, since every FIID is a weak* limit of
block FIIDs, it generally suffices to find an FIID with the desired
property on $\Td$ and to calculate $\Pbig{\phi(\cbuldot)(\Tdo) = 1}$. Indeed,
we may work with weak* limits of FIIDs.

We give two examples of this method that are inspired by \ref b.CGHV/.
They were the first to use Gaussian factors for similar purposes.

A \dfn{bisection} of a finite graph $\gh$ is a subset $S \subset \verts$
such that $\big| |S| - |\verts\setminus S| \big| \le 1$. In particular, if
$|\verts|$ is even, then $|S| = |\verts|/2$. The \dfn{size} 
of a bisection
$S$, written $\sz(S)$, is the number of edges $E(S, \verts \setminus S)$
that join $S$ to $\verts\setminus S$.
The problems of minimizing or maximizing the size of a bisection in a
regular graph are known to be hard in various senses and are of interest in
computer science; see \ref b.DDSW/. For a sequence of graphs $G_n$, define
$$
\minbi := \limsup_{n \to\infty} \min \big\{ \sz(S)/|\verts(G_n)| \st S \hbox{
is a bisection of } G_n \big\}
$$
and
$$
\maxbi := \liminf_{n \to\infty} \max \big\{ \sz(S)/|\verts(G_n)| \st S \hbox{
is a bisection of } G_n \big\}
\,.
$$
For random $d$-regular graphs, \ref b.Bollobas:iso/ proved that 
$$
\minbi \ge 
{d \over 4} - {\sqrt{d \log 2} \over 2}
=
{d \over 4} - 0.416^+ \sqrt d
\,,
$$
whereas \ref b.Alon:edge-exp/ proved that
$$
\minbi \le 
{d \over 4} - {3\sqrt{d} \over 32 \sqrt 2}
=
{d \over 4} - 0.0663^- \sqrt d
\,.
$$
We improve the latter (upper) bound to $d/4 - 0.32^-\sqrt d$.
Still in the context of random $d$-regular graphs,
\ref b.DMS:cuts/ establish the asymptotic values as $d \to\infty$
$$
\minbi
=
{d \over 4} - \parisi \sqrt{d \over 4} + o(\sqrt d)
$$
and
$$
\maxbi
=
{d \over 4} + \parisi \sqrt{d \over 4} + o(\sqrt d)
\,,
$$
where $\parisi \approx 0.7632$ is a certain known constant; for comparison
with the previous bounds, note that $\parisi/2 \approx 0.3816$.
The best previous results on $\minbi$ and $\maxbi$ for random $d$-regular
graphs for specific $d$
can be found in \ref b.MonPre/, \ref b.DDSW/, and \ref b.DSW/. In
the case of degrees $d = 3, 4$, we improve those results here, which were
that a.s., $\minbi \le 1/6$ and $\maxbi \ge 1.32595$ for $d = 3$ and
$\minbi \le 1/3$ and $\maxbi \ge 5/3$ for $d = 4$. We shall not actually
need that our finite graphs be regular.

\procl t.bisection
Let $G_n$ be finite graphs whose random weak limit is $\Td$ and whose
average degree tends to $d \ge 3$.
Then 
$$
\minbi < {d \over 2\pi} \arccos {2 \sqrt{d-1} \over d}
< {d \over 4} - {\sqrt d \over \pi}
\label e.gendmin
$$
and
$$
\maxbi > {d \over 2\pi} \arccos {-2 \sqrt{d-1} \over d}
> {d \over 4} + {\sqrt d \over \pi}
\,.
\label e.gendmax
$$
For $d = 3$, this gives
$$
\minbi < {3 \over 2\pi} \arccos \sqrt{8\over9} = 0.1622602^-
$$
and
$$
\maxbi > {3 \over 2\pi} \arccos \Bigl(-\sqrt{8\over9}\Bigr) = 1.3377398^+
\,;
$$
for $d = 4$, it yields
$\minbi < 1/3$ and $\maxbi > 5/3$.
\endprocl

\proof
We first make precise the connection of $\minbi$ and $\maxbi$ to weak*
limits of FIID processes on $\Td$. 
Let $\wf$ be the class of weak* limits $\phi$ of FIID processes on
$\Td$ with values $\phi(\cbuldot)(x) \in \{0, 1\}$ for $x \in \verts(\Td)$
and with $\P[\phi(\cbuldot)(\Tdo) = 1] =
1/2$.
Let $\minbifac$ be the infimum of
$\Ebig{ |\{x \sim \Tdo \st \phi(\cbuldot)(x) \ne \phi(\cbuldot)(\Tdo)\}|}/2$
taken over $\phi \in \wf$;
it is easily seen that this infimum is a minimum.
Similarly, 
let $\maxbifac$ be the supremum of $\Ebig{ |\{x \sim \Tdo \st \phi(\cbuldot)(x) \ne
\phi(\cbuldot)(\Tdo)\}|}/2$ taken over the same $\phi \in \wf$.
We claim that 
$$
\minbi \le \minbifac
\quad \hbox{ and } \quad
\maxbi \ge \maxbifac
\,.
\label e.emulate
$$

Indeed, let $\phi \in \wf$ and
$\epsilon > 0$.
There exists a block FIID $\phi_r \in \wf$ such that 
$$
\big|\Ebig{ |\{x \sim \Tdo \st \phi(\cbuldot)(x) \ne \phi(\cbuldot)(\Tdo)\}|}
-
\Ebig{ |\{x \sim \Tdo \st \phi_r(\cbuldot)(x) \ne \phi_r(\cbuldot)(\Tdo)\}|}\big|
<
\epsilon/2
\,.
$$
Now emulate $\phi_r$ on $G_n = (\verts_n, \edges_n)$.
Let $S_n \subseteq \verts_n$ be the subset of vertices assigned 1; this
need not be a bisection, as we know only that $\Ebig{|S_n|}/|\verts_n| \to
1/2$ as $n \to\infty$.
However, finite dependence of the block FIID $\phi_r$ implies that linear
deviations from the mean of $|S_n|$ are exponentially unlikely as $n
\to\infty$.
Furthermore, $\Ebig{\sz(S_n)}/|\verts_n|$ tends, as $n \to\infty$,
to $\Ebig{ |\{x \sim \Tdo \st \phi_r(\cbuldot)(x) \ne
\phi_r(\cbuldot)(\Tdo)\}|}/2$.
We have similar exponentially fast convergence for this 
proportion, $\sz(S_n)/|\verts_n|$.
Thus, for large $n$, there exists $S_n$ such that 
$$
\left|{|S_n| \over \verts_n} - {1 \over 2}\right|
< {\epsilon \over 8d}
$$
and 
$$
\left|{\sz(S_n) \over |\verts_n|} - 
{1\over2}\Ebig{ |\{x \sim \Tdo \st \phi_r(\cbuldot)(x) \ne \phi_r(\cbuldot)(\Tdo)\}|}\right|
< {\epsilon \over 4}
\,.
$$
In addition, $\lim_{n \to\infty} 2|\edges_n|/|\verts_n| = d$.

Now, if $|S_n| > |\verts(G_n) \setminus S_n| + 1$, remove the fewest number
needed of the
smallest-degree vertices in $S_n$ to obtain a bisection $S'_n$, while if
$|S_n| < |\verts(G_n) \setminus S_n| - 1$, add the
fewest number needed of the smallest-degree vertices
not in $S_n$ to obtain a bisection $S'_n$. The vertices moved from one part
to the other each have degree at most the median degree, which is at most
twice the mean degree, $4|\edges_n|/|\verts_n|$, whence for
large $n$, this new bisection $S'_n$
satisfies 
$$
\left|{\sz(S'_n) \over |\verts_n|}
- {1\over2} \Ebig{ |\{x \sim \Tdo \st \phi(\cbuldot)(x) \ne \phi(\cbuldot)(\Tdo)\}|}
  \right|
< \epsilon
\,,
$$
which proves \ref e.emulate/.

It remains to prove the asserted bounds but for $\minbifac$ and
$\maxbifac$.

Let $\sigma_{\pm}$ be the Gaussian wave functions of \ref b.CGHV/
with eigenvalues $\pm \rho_d$ for the transition operator.
These are weak* limits of FIID processes, and thus so are $\sgn \sigma_{\pm}$.
(The two wave functions are also related to each other via the distributional
equality $\sigma_+ \eqD f_* \sigma_-$, where $f \colon \R^\verts \to
\R^\verts$ is the map $\big(f(\omega)\big)(x) = (-1)^{d(\Tdo, x)}
\omega(x)$.)
Consider the bisections $\{x \in \verts \st \sigma_{\pm}(x) > 0\}$.
Now for jointly normal centered random variables $(Z_1, Z_2)$, we have
$$
\P[\sgn Z_1 \ne \sgn Z_2] = \P[Z_1 Z_2 < 0]
= {1 \over \pi} \arccos \Corr(Z_1, Z_2)
\,.
\label e.Gausspr
$$
Since $\Corr\big(\sigma_{\pm}(x), \sigma_{\pm}(y)\big) = \pm \rho_d$
for neighbors $x$ and $y$, 
we obtain that
$$
\minbifac \le {d \over 2\pi} \arccos \rho_d
$$
and
$$
\maxbifac \ge {d \over 2\pi} \arccos (-\rho_d)
\,.
$$

\efiginsnumber local-improve x5

Local improvements lead
to strict inequalities. That is, consider a vertex $x$ such as the heavily
circled one in \ref f.local-improve/, where the left figure applies to
$\sgn \sigma_+$ and the right figure to $\sgn \sigma_-$. Only the case of
$d=3$ is drawn, but all degrees are similar. It is easily checked that
such configurations have positive probability by using the Markov property
established in the proof of Theorem 3 of \ref b.CGHV/.
When such a configuration occurs, change the value at $x$ to its opposite;
likewise for configurations that are all opposite to those drawn. Note that
the lower neighbor of $x$ may change as well, but the upper neighbors of
$x$ will not. Thus, the number of edges incident to $x$ with the opposite
sign strictly decreases on the left and strictly increases on the right.

Finally, to prove the last inequalities in \ref e.gendmin/ and \ref
e.gendmax/ that involve an estimate of the arccos function, a little
algebra reveals that they are equivalent to the inequality 
$$
\sin {2 \over \sqrt d} < {2 \sqrt{d-1} \over d}
$$
for $d \ge 3$. Substituting $x := 2/\sqrt d$ shows that this is the same as
$\sin x < x \bigl(1 - x^2/4\bigr)^{1/2}$.
Indeed, $\sin x < x - x^3/6 + x^5/120$ for $0 < x^2 < 6$, whereas
$x \bigl(1 - x^2/4\bigr)^{1/2} > x(1 - x^2/7)$ for $0 < x^2 < 7/4$,
and this leads to the desired inequality.
\Qed

Finally, we improve \ref b.KKV/, who showed that if $G$ is a finite graph of
maximum degree 3 and girth at least 637,789, then there is a probability
measure on edge cuts of $G$ such that each edge belongs to a random cut
with probability at least 0.88672, whence (by taking expected size of edge
cuts) $G$ contains an edge cut of
cardinality at least $0.88672|\edges(G)|$.
For numerical comparison, note that this translates to the following result
when $G$ is 3-regular: 
3-regular $n$-vertex graphs of girth tending to infinity
possess subsets $S$ such that 
$\big|E(S,\verts\setminus S)\big| \ge \big(1.33008 - o(1)\big) n$ for even $n
\to\infty$.
\ref t.bisection/ already improved this by increasing the constant and by
requiring $S$ to be a bisection.

\procl t.edgecut
If $G$ is a finite graph of maximum degree $d$ and girth at least $2n+1$,
then there is a random edge cut $\cut$ of $G$ such that 
$$
\P[e \in \cut] \ge 
{1 \over \pi} \arccos {-\rho_d \over 1 + {d-1 \over d(n-1)}}
$$
for all $e \in \edges(G)$. 
\endprocl

For $d = 3$, this says, e.g., that if $G$ has girth at least 655, then
there is a random edge cut $\cut$ such that $\P[e \in \cut] \ge 0.89$ for
all $e\in \edges(G)$.

\proof
It suffices to prove the analogous result on $\Td$ via a block FIID of
radius $n$: We can then adjoin trees to $G$ in order to create a (possibly
infinite)
$d$-regular graph $G'$. The block FIID can be applied to $G'$ to obtain a
random cut $\cut'$ of $G'$; then we may let $\cut := \cut' \cap \edges(G)$.

To this end,
let $\Gm$ be standard Gaussian measure on $\R$. Then
the coordinate projections $Z_x \colon \big(\R^{\verts(\Td)},
\Gm^{\verts(\Td)}\big) \to \R$ are independent standard normal random
variables for $x\in \verts(\Td)$.
Put 
$$
F := {1 \over \sqrt{1+(n-1)d/(d-1)}} \sum_{d(\Tdo, x) < n}
{Z_x \over \big({-}\sqrt{d-1}\big)^{d(\Tdo, x)}}
\,.
$$
Then $F$ defines a block FIID $\sigma$ of radius $n$ with single marginal
equal to standard Gaussian.
Since 
$$
\Corr\big(\sigma(x), \sigma(y)\big)
=
-{1 \over 1+{(n-1)d \over d-1}} \sum_{k=0}^{n-2} {2 (d-1)^k \over
\sqrt{d-1}^{2k+1}}
=
-{\rho_d \over 1 + {d-1 \over d(n-1)}}
$$
for $x \sim y$,
the result follows by \ref e.Gausspr/. 
\Qed

\noindent
{\bf Acknowledgement.} I thank Allan Sly and Yuval Peres
for allowing me to include \ref t.sly/.
I also thank Lewis Bowen, Yuval Peres and Jeff Steif
for various discussions of this material.
I am grateful to a referee for helpful remarks.

\def\noop#1{\relax}
\input \jobname.bbl

\filbreak
\begingroup
\eightpoint\sc
\parindent=0pt\baselineskip=10pt

Department of Mathematics,
831 E 3rd St,
Indiana University,
Bloomington, IN 47405-7106 
\emailwww{rdlyons@indiana.edu}
{http://mypage.iu.edu/\string~rdlyons/}

\endgroup

\bye